\documentclass[12]{article}
\usepackage{amssymb,latexsym,amsmath}     
\usepackage[ngerman, english]{babel}
\usepackage[english]{babel}
\usepackage{latexsym,amssymb,colordvi}
\usepackage[utf8]{inputenc}
\usepackage{graphicx}
\date{}

\newtheorem{theorem}{Theorem}[section]

\newtheorem{qu}[theorem]{Question}

\title{\bf Note on the Number of Finite Groups of a Given Order}

\author{ \bf   A. R. Ashrafi\footnote{Corresponding author (Email: ashrafi@kashanu.ac.ir).} and E. Haghi   \\ Department of Pure Mathematics, Faculty of Mathematical Sciences,\\ University of Kashan, Kashan 87317$-$53153, I. R. Iran}

\begin{document}

\maketitle

\begin{abstract}
Let  $n$ be a positive integer and $G(n)$ denote the number of non-isomorphic finite groups of order $n$. It is well-known that $G(n) = 1$ if and only if $(n,\phi(n)) = 1$, where $\phi(n)$  and $(a, b)$ denote  the Euler's totient function and the greatest common divisor of $a$ and $b$, respectively. The aim of this paper is to first present a new proof for the case of $G(n) = 2$ and then give a  solution to the equation of $G(n) = 3$.

\vskip 3mm


\vskip 3mm

\end{abstract}

\section{Introduction}
Throughout this paper all groups are assumed to be finite. The function $G(n)$ and $\phi(n)$  denote the number of  non-isomorphic finite groups of order $n$ and the number of positive integers $i$ such that $1 \leq i \leq n$ and $(i,n) = 1$, respectively. For any other notation not defined here, we refer the reader to \cite{00,6}. Our calculations are done with the aid of GAP \cite{8}.

The problem of finding exact value of $G(n)$ in general seems to be very
difficult, but there are some very good estimations for this function in literature. M. R. Murty and V. K. Murty \cite{4},  used this well-known fact that groups of square-free order are supersolvable to prove that  for a square free number $n$,  $G(n) \leq \phi(n)$.

Dickson \cite{0} determined those positive integers $n$ for which every group of order $n$ is abelian. As a consequence of his result it can be proved that $G(n) = 1$ if and only if $(n,\phi(n)) = 1$. We encourage the interested readers to consult papers \cite{1,3} for a simpler proof  of Dickson's result. This result is our motivation to study the equation $G(n) = k$, for some small values of $k$. In an exact phrase,

\begin{qu}\label{q1}
For which positive integers $n$ does $G(n) = 2$? For which $n$ does $G(n) = 3$?
\end{qu}

This paper will be concerned with Question \ref{q1}. For the sake of completeness, we mention here some results which are crucial throughout this paper.

Given a group $G$ of square-free order greater than 1, a
theorem of H$\ddot{\rm o}$lder \cite{2} tells us that $G$ is metacyclic. If $n = p_1p_2\ldots p_s$, where $p_1 > p_2 > \ldots > p_s$ are primes then
\begin{equation}
G(n) = \sum_{S}\left(  \prod_{j=1}^r \frac{p_{S(j)}^{c_{S(j)}} - 1}{p_{S(j)} - 1} \right),
\end{equation}
where the summation is taken over  all the subsets $S = \{ S(1), 5(2), \ldots, S(r)\}$ of $\{ 2, 3,
. . . , n\}$ and $c_{S(j)}$ is the number of differences $p_i - 1, i \not\in S$, which are divisible by
$p$. If we define $f(n,m) = \prod_{q | m}(n,q-1)$, then one can see that $$G(n) = \sum_{d|n}\prod_{p | d}\frac{f(p,\frac{n}{d} - 1)}{p-1},$$ where $p$ and $q$ are primes, see \cite{000} for details.

Suppose $G$ is a finite group and $H$ is a normal subgroup such that $p \nmid |H|$ and $\frac{|G|}{|H|} = p^\alpha$, where $\alpha$ is a positive integer. Then $H$ is called a normal $p-$complement for $G$. In \cite[p. 252]{15}, the following result is proved:

\begin{theorem}
(Burnside's Transfer Theorem). Let $G$ be a finite group and $P$ be a Sylow $p-$subgroup of $G$. If $P \leq Z(N_G(P))$ then $G$ has a normal $p-$complement.
\end{theorem}

Suppose $n = \prod_{i=1}^rp_i^{\alpha_i}$ and define $\Phi(n) =\prod _{i=1}^r(p^{\alpha_i} - 1)$. Redei \cite[Satz 10]{7}, characterized the  positive integers $n$ with this property that all abelian groups of order $n$ are abelian. Such a positive integer is called an abelian number. He proved that:

\begin{theorem}\label{thm3}
 $n = \prod_{i=1}^rp_i^{\alpha_i}$ is an abelian number if and only if $(n,\Phi(n)) = 1$ , for each $i$, $1 \leq i \leq r$.
\end{theorem}

In the same manner,  we call a positive integer $n$ a nilpotent number if every group of order $n$ is nilpotent. Pakianathan and Shankar \cite{5}, proved that a positive integer $n = p_1^{\alpha_1}p_2^{\alpha_2} \cdots p_3^{\alpha_t}$, $p_i$'s are distinct primes, is a nilpotent number if and only if $p_i^k \not\equiv 1 \ (mod p_j)$, for all integers $i,$ $j$ and $k$ with $1 \leq k \leq \alpha_i$.

In this paper, our interest is the study of the function $G(n)$. The main result of this paper is as follows:

\begin{theorem} \label{thm4} Suppose $n = \prod_{i=1}^rp_i^{\alpha_i}$.
\begin{enumerate}

\item $G(n) = 2$ if and only if one of the following are hold:
\begin{enumerate}
\item $r = 1$ and $\alpha_1 = 2$;

\item $r \geq 2$, $n$ is square-free and there is a unique pair $(i,j)$ such that $p_i | p_j -1$;

\item $n = p_1^2p_2\ldots p_r$, $r \geq 2$, and for each $i, j$ such that $1 \leq i,j \leq k$, we have $p_i \nmid p_j - 1$; and for each $k$, $2 \leq k \leq r$, $p_k \nmid p_1^2 - 1$.
\end{enumerate}

\item $G(n) = 3$ if and only if one of the following are hold:
\begin{enumerate}
\item $r \geq 3$, $n$ is square-free and there is a unique triple $(i,j,l)$ such that $p_i | p_j -1$ and $p_j$ $|$ $p_l - 1$;

\item $n = p_1^2p_2\ldots p_r$, $r \geq 2$, and for each $i, j$ such that $1 \leq i,j \leq r$, we have $p_i \nmid p_j - 1$; and there exists a unique positive integer $l$, $2 \leq l \leq r$ such that $p_l | p_1^2 - 1$.
\end{enumerate}

\end{enumerate}
\end{theorem}

\section{Proof of Main Result}
The aim of this section is to prove our main results. At first, we note  that there are exactly five groups of order $p^3$, where $p$ is prime. It is well-known that $Z_8$, $Z_2 \times Z_4$, $Z_2 \times Z_2 \times Z_2$, $D_8$ and $Q_8$ are all groups of order eight and if $p$ is an odd prime then $Z_{p^3}$, $Z_p \times Z_{p^2}$ and $Z_p \times Z_p \times Z_p$ are only abelian groups of order $p^3$ and there are exactly two  non-abelian groups of this order as follows:
\begin{enumerate}
\item $G_1 = \langle x, y \ | \ x^p = y^{p^2} = 1, x^{-1}yx = y^{p+1}     \rangle$,

\item $G_2 = \langle x, y, z \ | \ x^p = y^{p} = z^p = 1, yz = zy, zx = xz, x^{-1}yx = yz    \rangle$.
\end{enumerate}
This proves that $G(p^3) = 5$. The number of groups of order $p^2q$, $p$ and $q$ are primes, are recorded in Table 1. This table is first appeared in the famous book of Burnside \cite{00}.

\vskip 3mm

\centerline{\textbf{Table 1.} The Number of Groups of Order $p^2q$, $p$ and $q$ are primes. }
\begin{center}
\begin{tabular}{l|l}
Conditions & $G(p^2q)$\\ \hline
$p \nmid q - 1$ $\&$ $q \nmid p - 1$ & 2\\
$p | q - 1$, $p^2 \nmid q - 1$ and $q \nmid p^2 - 1$ & 4\\
$p^2 | q - 1$ & 5\\
$q = 2$ and $q | p - 1$ & 5\\
$q$ is odd and $q | p - 1$ & $\frac{q+1}{2} + 4$
\end{tabular}
\end{center}

Suppose $p$ is prime and $q | p - 1$. The Frobenius group $F_{p,q}$ can be presented as follows:
$$F_{p,q} = \langle a, b \ | \ a^p = b^q = 1, b^{-1}ab = a^u\rangle,$$
where $u$ is an element of order $q$ in the unit group of the ring $Z_p$, \cite[p. 290]{24}.

\subsection{The Case that $G(n) = 2$.}
Suppose $G(n) = 2$, where $ n=\prod_{i=1}^r p_i^{r_i} $. If there exists $i$ such that $ r_i\geq3 $ then there are five non-isomorphic finite groups of order $p_i^3$. This shows that we have at least five groups
$  Z_{p_i^{r_i-3}p_1^{r_1}\ldots p_{i-1}^{r_{i-1}}p_{i+1}^{r_{i+1}}\ldots p_{k}^{r_{k}}} \times A$, where $A$ is a group of order $p_i^3$. Also, if there are $i$ and $j$, $i \ne j$ and $r_i, r_j \geq 2$ then $$Z_{\frac{n}{p_i^2p_j^2}} \times Z_{p_i}\times  Z_{p_i}\times  Z_{p_j^2},
Z_{\frac{n}{p_i^2p_j^2}} \times Z_{p_i^2p_j^2}, Z_{\frac{n}{p_i^2p_j^2}} \times  Z_{p_i^2}\times  Z_{p_j}\times  Z_{p_j}, Z_{\frac{n}{p_i^2p_j^2}} \times  Z_{p_i}\times  Z_{p_i}\times  Z_{p_j}\times  Z_{p_j}$$ are four non-isomorphic finite groups of order
$n$ which is not possible. Thus, $n$ is square-free or $n$ is a cube-free number such that there exists exactly one prime number $p$ such that $p^2 | n$. If $n$ is square-free number then by Equation (1), there exists only one pair $(i,j)$ such that $p_i | p_j - 1$, as desired. So, we can assume that $n$ is a cube-free number such that there exists exactly one prime number $p$ such that $p^2 | n$. Since $Z_n$ and $Z_{p} \times Z_{\frac{n}{p}}$ are two non-isomorphic abelian group of order $n$, every group of order $n$ has to be abelian. By Theorem \ref{thm3}, the condition $1(c)$ in Theorem \ref{thm4} is satisfied.

Note that if there are only two finite groups of order $n$ then $(n,\phi(n))$ is primes, but the converse is not generally correct. To see this, it is enough to check that $G(75) = 3$ and $(75,\phi(75)) = (75,40) = 5$.

\subsection{The Case that $G(n) = 3$.}
In this case the finite groups with $G(n) = 3$ will be characterized. By the proof of Theorem \ref{thm4}(1), such a number $n$ is square-free or it is a cube-free integer that there exists a unique prime number $p$ such that $p^2 | n$. Our main proof will consider two different cases as follows:
\begin{enumerate}
\item $n = p_1p_2 \cdots p_r$, $r \geq 3$, is square-free. By Equation (1), there are exactly two pairs $(i,j)$ and $(i^\prime,j^\prime)$ such that $p_i | p_j - 1$ and  $p_{i^\prime} | p_{j^\prime} - 1$.
\begin{enumerate}
\item \textit{$p_{i} = p_{i^\prime}$ and $p_{j} \ne p_{j^\prime}$}. In this case, by Equation (1), there are exactly $p_i + 2 \geq 4$ finite groups of order $n$, which is a contradiction.

\item \textit{$p_{i} \ne p_{i^\prime}$ and $p_{j} = p_{j^\prime}$}. Under this conditions there are four groups $Z_n$, $Z_{\frac{n}{p_ip_j}}$ $\times$ $F_{p_j,p_i}$, $Z_{\frac{n}{p_i^\prime p_j^\prime}}$ $\times$ $F_{p_j^\prime,p_i^\prime}$ and $Z_{\frac{n}{p_ip_jp_i^\prime}}$ $\times$  $F_{p_j,p_ip_i^\prime}$ are four non-isomorphic groups of order $n$, which is impossible.

\item $p_{i} = p_{j^\prime}$. By Equation (1), there are exactly three finite groups of order $n$, as desired.
\end{enumerate}

\item \textit{$n=p_1^2p_2\ldots p_r $, where $p_i$'s are distinct prime integers and $r \geq 2$.} Suppose there are $i$ and $j$, $2 \leq i, j \leq r$, such that $p_i | p_j - 1$. Then there are four non-isomorphic groups as follows:
\begin{eqnarray*}
    &&Z_n,\\ &&Z_{\frac{n}{p_1^2p_ip_j}} \times Z_{p_1^2} \times  F_{p_j,p_i},\\ &&Z_{\frac{n}{p_1^2p_ip_j}} \times Z_{p_1} \times Z_{p_1} \times  Z_{p_ip_j}\\ &&Z_{\frac{n}{p_1^2p_ip_j}} \times Z_{p_1} \times Z_{p_1} \times  F_{p_j,p_i},
\end{eqnarray*}
which is not possible. If there exists $i$, $2 \leq i \leq r$, such that $p_1 | p_i - 1$ or $p_i | p_1 - 1$ then by Table 1, there are at least four groups of order $n$ leads again to a contradiction. We now assume that the prime integers $p_1, \cdots,$ $p_r$ are satisfied the condition $2(b)$ in Theorem \ref{thm4}. Choose $P_l$  to be the Sylow $p_{_l}-$subgroup of $G$ and $N_l = N_G(P_l)$. Since $(|\frac{N_l}{P_l}|,|Aut(P_l)|) = 1$, $P_l \leq Z(N_l)$ and by Burnside'e transfer theorem $G$ has a normal subgroup $N$ of order $p_1^2p_2\ldots p_{l-1}p_{l+1}\ldots p_k$ = $p_1s$. By Part (1) of this theorem, $N \cong Z_{p_1s}$ or $Z_{p_1} \times Z_{s}$. Since
\begin{eqnarray*}
Aut(Z_{p_1s}) &=& Z_{p_1^2-p_1}\times Z_{p_2-1}\times \ldots \times   Z_{p_{l-1}-1}\times  Z_{p_{l+1}-1}\times \ldots \times  Z_{p_k-1},\\
Aut(Z_{p_1}\times Z_s)&=& GL(2,p_1)\times Z_{p_2-1}\times \ldots \times  Z_{p_{l-1}-1}\times  Z_{p_{l+1}-1}\times \ldots \times  Z_{p_k-1}
\end{eqnarray*}
and $p_l \nmid | Aut(Z_{p_1s}) |  $, there exists a unique homomorphism $\varphi_1$ from $P_l$ into $Aut(Z_{p_1s})$. This homomorphism has to be trivial. On the other hand, $ |GL(2,p_1)|=(p_1^2-p_1)(p_1^2-1) $,
$ p_l \mid  | Aut(Z_{p_1}\times Z_s) |  $ and all Sylow $p_l$ subgroups of $GL(2,p_1)$ are cyclic. This shows that all subgroups of order $p_l$ are conjugate in $GL(2,p_1)$. If $\varphi$ and $\psi$ are two non-trivial homomorphism from $P_l$ into $Aut(Z_{P_1} \times Z_s)$ then
$ P_l\times _ \psi (Z_{p_1}\times Z_s)\cong P_l\times _ \varphi (Z_{p_1}\times Z_s) $. Therefore, there are three groups of order $n$ which can be presented as semi-direct product by the following homomorphism:
\begin{eqnarray*}
\varphi_1:P_l & \rightarrow & {Aut(Z_{p_1s})}\\
x&\rightarrow & 1\\ [.3cm]
\varphi_2:P_l & \rightarrow & {Aut(Z_{p_1}\times Z_s)}\\
x&\rightarrow & 1\\ [.3cm]
\varphi_3:P_l &\rightarrow & {Aut(Z_{p_1}\times Z_s)} \\
x&\rightarrow & y
\end{eqnarray*}
Here, $x$ and $y$ are two elements of order $p_l$ in the groups $P_l$ and $Aut(Z_{p_1}\times Z_s)$, respectively.

\end{enumerate}
This completes the proof.

\section{Concluding Remarks}

In this paper the equation $G(n) = 2, 3$ is considered into account. In the case that $G(n) = 2$, we present a new proof for an old result of M. R. Murty and V. K. Murty. Then we continue the pioneering work of these authors to solve the equation $G(n) = 3$. We check also the case for $G(n) = 4$, but the problem for this case is so difficult. So, the complete solution of $G(n) = 4$ is a problem for future.

\vskip 3mm

\noindent{\bf Acknowledgement.} The research of the authors
are partially supported by the University of Kashan under grant no
364988/101.

\end{document}